\documentclass[number,citesort,MSNbibl,seceqn,dvips]{arxbj}
\usepackage{multirow}

\aid{0}
\volume{18}
\issue{2}
\pubyear{2012}
\firstpage{586}
\lastpage{605}
\doi{10.3150/10-BEJ342}

\makeatletter
\newremark{example}{Example}
\newremark{rem}{Remark}[section]
\newtheorem{lem}{Lemma}[section]
\newtheorem{theorem}{Theorem}[section]
\renewcommand{\citet}{\cite}
\makeatother

\begin{document}
\begin{frontmatter}

\title{Reparametrization of the least favorable submodel in semi-parametric multisample~models}
\runtitle{Reparametrization of the least favorable submodel}

\begin{aug}
\author[1]{\fnms{Yuichi} \snm{Hirose}\corref{}\thanksref{1}\ead[label=e1]{Yuichi.Hirose@msor.vuw.ac.nz}}
\and
\author[2]{\fnms{Alan} \snm{Lee}\thanksref{2}\ead[label=e2]{lee@stat.auckland.ac.nz}}
\runauthor{Y. Hirose and A. Lee}
\address[1]{School of Mathematics, Statistics and Operations Research, Victoria University of Wellington,
 New Zealand. \printead{e1}}
\address[2]{Department of Statistics, University of Auckland, New~Zealand. \printead{e2}}
\end{aug}

\received{\smonth{4} \syear{2010}}
\revised{\smonth{10} \syear{2010}}

\begin{abstract}
The method of estimation in Scott and Wild (\textit{Biometrika}
\textbf{84} (1997) 57--71 and \textit{J. Statist. Plann. Inference}
\textbf{96} (2001) 3--27) uses a reparametrization of the profile
likelihood that often reduces the computation times dramatically.
Showing the efficiency of estimators for this method has been a
challenging problem.
In this paper, we try to solve the problem by investigating conditions
under which the efficient score function and the efficient information
matrix can be expressed in terms of the parameters in the
reparametrized model.
\end{abstract}


\end{frontmatter}

\section{Introduction}\label{sec1}
In a series of papers, Scott and Wild \citet{ScottWild1997,ScottWild2001} developed
methods of reparametrization of profile likelihood that can be applied
to a variety of response-selective sampling designs.
The advantage of the methods is that they often give us computationally efficient estimators. 
The (statistical) efficiency of these methods has been demonstrated in
special cases by several authors. For example,
Breslow, Robins and Wellner \citet{BreslowRobinsWellner2000} considered case-control sampling where
either a case or control is selected by a randomization device with
known selection probabilities, and the covariates of the resulting case
or control are measured. In the case of two-phase, outcome-dependent
sampling, Breslow, McNeney and Wellner \citet{BreslowMcNeneyWellner2003} applied the missing value
theory of Robins, Rotnitzky and Zhao \citet{RobinsRotnitzkyZhao1994}
and Robins, Hsieh and Newey \citet{RobinsHsiehNewey1995}. Here, individuals in the population are
selected at random and their status (e.g., case or control) is
determined. Then, with a probability depending on their status, the
covariates are measured. The unobserved covariates are treated as
missing data. Lee and Hirose \citet{LeeHirose2008} used the profile likelihood method
to derive a semi-parametric efficiency bound, and then showed that this
bound coincides with the asymptotic variance of the Scott--Wild
estimator, hence demonstrating the efficiency of the estimator.

In Lee and Hirose \citet{LeeHirose2008}, it was demonstrated that, in the case of the
Scott--Wild estimator, it is possible to reparametrize the least
favorable submodel so that the efficient score function and the
efficient information matrix can be expressed in terms of the
parameters in the reparametrized model.

The aim of this paper is to investigate conditions under which a
reparametrization of the least favorable submodel yields an efficient
estimation.

We consider an $S$-vector of semi-parametric models $(\mathcal{P}_1,\ldots,\mathcal{P}_S)$ where, for each $s=1,\ldots,S$,
\[
\mathcal{P}_s=\{p_s(x;\beta,\eta)\dvt\beta \in \Theta_{\beta}\subset R^m, \eta \in  \Theta_{\eta}\}
\]
is a probability model on the sample space $\mathcal{X}_s$ with the
parameter of interest $\beta$, an $m$-dimensional parameter, and the
nuisance parameter $\eta$, which may be an infinite-dimensional
parameter. Let $(\beta_0,\eta_0)$ be the true value of $(\beta,\eta)$.
We assume $\Theta_{\beta}$ is a~compact set containing an open
neighborhood of $\beta_0$ in $R^m$, and $\Theta_{\eta}$ is a convex set
containing~$\eta_0$ in a Banach space $\mathcal{B}$. We refer to the
$S$-vector of semi-parametric models $(\mathcal{P}_1,\ldots,\mathcal{P}_S)$
as the multisample model.

Under the model, we observe $S$ independent samples $X_{s1},\ldots,X_{s
n_s}$ ($s=1,\ldots,S$), where
 $X_{s1},\ldots,X_{s n_s}$ are independently and identically distributed (i.i.d.) according to the model $\mathcal{P}_s$.
Let $n=\sum_{s=1}^S n_s$. We assume the sample size proportions
$(n_1/n,\ldots,n_S/n)$ converge to weight probabilities
$(w_1,\ldots,w_S)$:
%
\begin{equation}\label{EqnSamplesizeProp}
\biggl(\frac{n_1}{n},\ldots,\frac{n_S}{n}\biggr) \rightarrow
(w_1,\ldots,w_S),
\end{equation}
where  $w_s>0$ and $\sum_{s=1}^S w_s =1$.

The log-likelihood for the multisample data is
\begin{equation}\label{log-likelihood_multi}
\ell_n(\beta,\eta)=\sum_{s=1}^S \sum_{i=1}^{n_s} \log
p_s(X_{si};\beta,\eta).
\end{equation}

The paper is organized as follows: In the rest of Section \ref{sec1}, we
give examples of semi-parametric multisample models. In Section
\ref{sec2}, we introduce the least favorable submodel in multisample models and
in Section \ref{Defreparameterization}, we present the main result of conditions under
which reparametrization gives efficient estimators in  multisample
models. In Section \ref{sec4}, we give a numerical example
and use the result developed in the paper  to show that the estimators
in the example are efficient.

\subsection{Examples}
The idea of multisample data is familiar from elementary statistics;
for example, the well-known two-sample $t$-test and the one-way ANOVA
for comparing several means both involve multiple samples. Following
are several semi-parametric examples.

\begin{example}[(Biased sampling model)]\label{example1} Vardi \cite{Vardi1985}
developed the method of estimation in the $S$-sample biased sampling
model with known selection bias weight functions. The following setup
and notation are from \cite{GillVardiWellner1988}.

Suppose that non-negative weight functions $w_1(x),\ldots,w_S(x)$ are
given and let $G(x)$ be an unknown distribution function on a sample
space $\mathcal{X}$. Define the corresponding biased sampling model by
\[
p_s(x;G)=\frac{w_s(x)g(x)}{W_s(G)} \qquad (s=1,\ldots,S),
\]
where $g(x)=\mathrm{d}G(x)/\mathrm{d}\mu$ with respect to Lebesgue measure $\mu$ and
$W_s(G)=\int_\mathcal{X} w_s(x)\,\mathrm{d}G(x)$. The $S$-sample biased sampling
model generates $S$ independent samples
\[
X_{s1},\ldots,X_{sn_s} \sim p_s(x;G)\qquad   (s=1,\ldots,S).
\]
%

Gilbert, Lele and Vardi \cite{GilbertLeleVardi1999} considered an extension of this model that
allows the weight function to depend on an unknown finite-dimensional
parameter $\theta$.

Suppose a set of non-negative weight functions
$w_1(x,\theta),\ldots,w_S(x,\theta)$ depend on $\theta$. The
semi-parametric biased sampling model is defined by
$$
p_s(x;\theta,G)=\frac{w_s(x,\theta)g(x)}{W_S(\theta,G)}\qquad  (s=1,\ldots,S),
$$
where $W_s(\theta,G)=\int_\mathcal{X}w_s(x,\theta)\,\mathrm{d}G(x)$.
Gilbert \cite{Gilbert2000} provides a large sample theory of this example.
\end{example}

The following examples are semi-parametric multisample models that all
have the same
   underlying data-generating process on the sample space $\mathcal{Y} \times \mathcal{X}$, called the full data model,
\[
\mathcal{Q}=\{p(y,x;\theta,G)=f(y|x;\theta)g(x)\dvt
\theta \in \Theta, G \in \mathcal{G} \},
\]
where $f(y|x;\theta)$ is a conditional density of $Y$ given $X$ that
depends on  a finite dimensional parameter $\theta$ and $G(x)$ is an
unspecified distribution function of $X$ that is an
infinite-dimensional nuisance parameter ($g(x)$ is the density of
$G(x)$). We assume the set $\Theta$ is a compact set containing a
neighborhood of the true value $\theta_0$ and $\mathcal{G}$ is the set of
all distribution functions of $x$.
Unless stated otherwise, $Y$ may be a discrete or continuous variable.

\begin{example}[(Case-control study)]\label{example2}
We assume that $Y$ takes
values in $\{1,\ldots,S\}$. In a~case-control study, due to the design,
we do not observe a random sample from the full data model $\mathcal{Q}$.
Instead, for each $s=1,\ldots,S$, we observe $n_s$-samples from the
conditional distribution $P(X|Y=s)$. By Bayes' theorem, the density of
$P(X|Y=s)$ is
\[
\frac{f(s|x;\theta)g(x)}{\int f(s|x;\theta)\,\mathrm{d}G(x)}.
\]
The case-control study is a special case of the semi-parametric biased
sampling model of Example \ref{example1} with weight functions
$w_s(x,\theta)=f(s|x;\theta)$ ($s=1,\ldots,S$).\vadjust{\goodbreak}
\end{example}

\begin{example}[(Missing data)]\label{example3}
Instead of observing full data
$(Y,X)$ from the full data model $\mathcal{Q}$ for all individuals, we
observe $(Y,X)$ for $n_0$-samples and observe $Y$ for $n_1$-samples.
The result is the multisample data
\[
(x_{01},y_{01}),\ldots,(x_{0n_0},y_{0n_0}), y_{11},\ldots,y_{1n_1}
\]
from a multisample model with densities
\[
p_0(y,x;\theta,g)=f(y|x;\theta)g(x)
\]
and
\[
p_1(y;\theta,g)=\int f(y|x;\theta)g(x)\,\mathrm{d}x.
\]
This example is not a special case of Example \ref{example1}.
\end{example}

\begin{example}[(Standard stratified sampling and two-phase, outcome-dependent sampling)]\label{example4}
For\vspace*{2pt} a~partition of the sample space $\mathcal{Y} \times \mathcal{X}=\bigcup_{s=1}^S \mathcal{S}_s$, let
\[
Q_s(\theta,G)=\int f(y|x;\theta) 1_{(y,x) \in \mathcal{S}_s}\,\mathrm{d}y\,\mathrm{d}G(x)
\]
be the probability of $(Y,X)$ belonging to stratum $\mathcal{S}_s$.

In standard stratified sampling, for each $s=1,\ldots,S$, a random
sample of size $n_s$ is taken from the conditional distribution
\[
p_s(y,x;\theta,G)=\frac{f(y|x;\theta)g(x)1_{(y,x)\in \mathcal{S}_s}}{Q_s(\theta,G)}
\]
of $(Y,X)$ given stratum $\mathcal{S}_s$.
This is a more general version of the semi-parametric biased sampling
model of Example \ref{example1} with weight functions
$w_s(y,x,\theta)=f(y|x;\theta)1_{(y,x)\in \mathcal{S}_s}$
($s=1,\ldots,S$).

Lawless, Kalbfleisch and Wild \cite{LKW1999} discussed variations of the two-phase, outcome-dependent
sampling design (the variable probability sampling designs (VPS1, VPS2)
and the basic stratified sampling design (BSS)).
For all sampling schemes (VPS1, VPS2 and BSS), we have $m_s$ fully
observed units and $n_s-m_s$ subjects where the only information
retained is the identity of the stratum, $s=1,\ldots,S$. The
corresponding likelihood is
%
\begin{equation}\label{EqnLikelihood}
L(\theta,G)= \Biggl\{ \prod_{s=1}^S \prod_{i=1}^{m_s}
f(y_{si}|x_{si};\theta)g(x_{si})\Biggr\}\Biggl\{
\prod_{s=1}^SQ_{s}(\theta,G)^{n_s-m_s}\Biggr\}.
\end{equation}
We interpret the observed data from two-phase, outcome-dependent
sampling as data from a multisample model with densities
\[
p_1(y,x;\theta,G)=f(y|x;\theta)g(x)
\]
and
\[
p_2(s;\theta,G)=Q_{s}(\theta,G).
\]
This example is not a special case of Example \ref{example1}.\vspace*{-2pt}
\end{example}

\section{The least favorable submodel}\label{sec2}\vspace*{-2pt}


The \textit{log-likelihood function} for a single observation in the
multisample model is
%
\begin{equation}\label{EqnLogLikelihoodOne}
\ell(s,x;\beta,\eta)= \log p_s(x;\beta,\eta)\qquad   (x \in \mathcal{X}_s,\
s=1,\ldots,S).
\end{equation}

The expectation with respect to the density $p_s(x;\beta,\eta)$ is
denoted by $E_{s,\beta,\eta}$.

We assume that there is a differentiable function $\beta \rightarrow
\hat{\eta}_{\beta}$ such that
%
\begin{equation}\label{ConditionR1}
\hat{\eta}_{\beta_0}=\eta_0
\end{equation}
and
%
\begin{equation}\label{ConditionR2}
\dot{\ell}^*(s,x)=\frac{\partial }{\partial
\beta}\bigg|_{\beta=\beta_0} \ell(s,x,\beta; \hat{\eta}_{\beta})
\end{equation}
is the efficient score function (definition of the efficient score
function in the multisample model is given in Appendix \ref{appendixa}). We call the
model
\[
p_s(x;\beta,\hat{\eta}_{\beta})\qquad  (\beta \in \Theta_{\beta},  s=1,\ldots,S),
\]
the least favorable submodel for the multisample model $(\mathcal{P}_1,\ldots,\mathcal{P}_S)$.\vspace*{-2pt}

\begin{rem}
Under mild regularity conditions with the assumption that
\[
\hat{\eta}_{\beta}=\operatorname{arg\,max}\limits_{\eta \in \Theta_{\eta}}\sum_{s=1}^S
w_s E_{s,\beta_0,\eta_0}\{\log p_s(X;\beta,\eta)\}
\]
exists for all $\beta$ in some neighborhood of $\beta_0$,
 (\ref{ConditionR2}) is the efficient score function due to~\cite{Newey1994}.
The definition of the least favorable submodel given above includes
this as a special case but we do not limit our consideration only in
this case.\vspace*{-2pt}
\end{rem}

Our approach uses the method in Scott and Wild \citet{ScottWild1997,ScottWild2001} to
find a candidate function~$\hat{\eta}_{\beta}$ as well as Theorem
\ref{theorema2}
in Appendix \ref{appendixa} to verify that (\ref{ConditionR2}) with the candidate
function gives the efficient score function. In the next example we
illustrate this procedure.\vspace*{-2pt}

\subsection{Example: Stratified sampling (continued)}\vspace*{-2pt}

Stratified sampling was introduced in Example \ref{example4}.

Let
\[
Q_{s|X}(x;\theta)=\int f(y|x;\theta) 1_{(y,x) \in
\mathcal{S}_s}\,\mathrm{d}y.\vadjust{\goodbreak}
\]

For each $s=1,\ldots,S$, let $F_{s0}$ be the cumulative distribution
function for the density $p_s(y,x;\theta_0,g_0)$ at the true value
$(\theta_0,g_0)$. The expected likelihood in the model is
\[
\sum_{s=1}^{S}w_s E_{s,0} \{\log p_s(y,x;\theta,g)\}
=\sum_{s=1}^{S}w_s \int \log p_s(y,x;\theta,g)\,\mathrm{d}F_{s0}(y,x) .
\]


For each $\theta$, the method in Scott and Wild \citet{ScottWild1997,ScottWild2001}
finds a maximizer $\hat{g}_{\theta}(x)$ of log-likelihood under the
assumption that the support of the distribution of $X$ is finite; that
is, $\textsc{supp}(X)=\{v_1,\ldots,v_K\}$. Let
$(g_1,\ldots,g_K)=\{g(v_1),\ldots,g(v_K)\}$. Then $\log g(x)$ and
$Q_s(\theta,g)$ can be expressed\vspace*{-1pt} as $\log g(x)=\sum_{k=1}^K 1_{x=v_k}
\log g_k  $ and $Q_s(\theta,g)=\int Q_{s|X}(x;\theta) g(x)\,\mathrm{d}x =
\sum_{k=1}^K Q_{s|X}(v_k;\theta)g_k$.\vspace*{1pt}

 To find the maximizer $(g_{1},\ldots,g_{K})$ of the expected log-likelihood
\[
\sum_{s=1}^{S}w_s \int \log p_s(y,x;\theta,g)\,\mathrm{d}F_{s0} =
\sum_{s=1}^S w_s \biggl[  \int  \{\log f(y|x;\theta) +  \log
g(x)\}\,\mathrm{d}F_{s0}  -  \log Q_s(\theta,g)\biggr]
\]
at $\theta$, differentiate this expression with respect to $g_k$ and
set the derivative equal to zero,
\[
\frac{\partial}{\partial  g_k  } \sum_{s=1}^{S}w_s \int \log p_s(y,x;\theta,g)\,\mathrm{d}F_{s0}
=\sum_{s=1}^S w_s \biggl\{ \frac{\int  1_{x=v_k}\,\mathrm{d}F_{s0} }{g_k} -
\frac{Q_{s|X}(v_k;\theta)}{Q_s(\theta,g)} \biggr\}=0.
\]
The solution $g_k$ to the equation is
\[
 \hat{g}_{\theta}(v_k)= g_k
  =  \frac{\sum_{s=1}^S w_s   \int  1_{x=v_k}\,\mathrm{d}F_{s0}}
{\sum_{s=1}^S w_s Q_{s|X}(v_k;\theta)/Q_s(\theta,g)}.
\]
%
The form of the function motivates us to prove the following result.

\begin{lem}[(The least favorable submodel)]\label{lem21}
For $\theta \in \Theta$, let
%
\begin{equation}\label{EqnhatgstratifiedSampling}
\hat{g}_{\theta}(x)
  =  \frac{f^*_0(x)}
{ \sum_{s=1}^S  w_s  Q_{s|X}(x;\theta)/\hat{Q}_{s}(\theta)},
\end{equation}
where
%
\begin{eqnarray}\label{Eqnf_star}
f^*_0(x)
=\sum_{s=1}^S  w_s
\frac{Q_{s|X}(x;\theta_0)g_0(x)}{Q_{s}(\theta_0,g_0)},
\end{eqnarray}
%
and
%
\begin{equation}\label{EqnhatQstratifiedSampling}
\hat{Q}_s(\theta)=
  \int Q_{s|X}(x;\theta) \hat{g}_{\theta}(x)\,\mathrm{d}x\qquad (s=1,\ldots,S).
\end{equation}
Then the efficient score function is given by
%
\begin{equation}\label{EfficientScore}
\dot{\ell}^*(s,y,x)=\frac{\partial}{\partial
\theta}\bigg|_{\theta=\theta_0} \log p_s(y,x;\theta,\hat{g}_{\theta}).
\end{equation}
\end{lem}

\begin{pf}
In Appendix \ref{appendixb}, we show that  $\sum_{s=1}^{S}w_s \int \log
p_s(y,x;\theta,\hat{g}_{\theta})\,\mathrm{d}F_{s0}$ satisfies
conditions~(\ref{eqnEffCond1}) and (\ref{eqnEffCond2}) in Theorem \ref{theorema2} in Appendix
\ref{appendixa} so that the claim follows from this theorem.\vspace*{-2pt}
\end{pf}

\begin{rem}\label{rem22}
 Note that equations (\ref{EqnhatgstratifiedSampling}) and (\ref{EqnhatQstratifiedSampling}) are consistent at $\theta=\theta_0$:
(\ref{EqnhatgstratifiedSampling}) and (\ref{Eqnf_star}) imply that
$\hat{g}_{\theta_0}(x)=g_0(x)$ if
$\hat{Q}_s(\theta_0)=Q_s(\theta_0,g_0)$. On the other hand, if
$\hat{g}_{\theta_0}(x)=g_0(x)$, we have $\hat{Q}_s(\theta_0)= \int
Q_{s|X}(x;\theta_0) g_0(x)\,\mathrm{d}x=Q_s(\theta_0,g_0)$ by
(\ref{EqnhatQstratifiedSampling}).\vspace*{-2pt}
\end{rem}

\section{Main result}\vspace*{-2pt}
\label{Defreparameterization}

Suppose there is a finite-dimensional, vector-valued function $\beta
\rightarrow q_{\beta}$ such that the density for the least favorable
submodel is of the form
%
\begin{equation}\label{EqnDefreparameterized1}
 p_s(x;\beta,\hat{\eta}_{\beta})= p^*_s(x;\beta,q_{\beta})\qquad
 \mbox{for all } \beta \in \Theta_{\beta}\ (s=1,\ldots,S),
\end{equation}
where the function $p^*_s(x;\beta,q)$ is twice continuously
differentiable with respect to $(\beta,q)$ and $q$ is a
finite-dimensional parameter. Further, suppose
%
\begin{equation}\label{EqnDefreparameterized2}
\sum_{s=1}^S w_s \int p^*_s(x;\beta,q)\,\mathrm{d}x=1\qquad  \mbox{for all }
(\beta,q) \in \Theta_{\beta} \times D_q,
\end{equation}
where $\Theta_{\beta}$ and $D_q$ are neighborhoods of $\beta_0$ and
$q_{\beta_0}$, respectively.
Then the model 
\[
p^*_s(x;\beta,q)\qquad(\beta \in \Theta_{\beta}, q \in D_q,s=1,\ldots,S),
\]
is called a \textit{reparametrized model} for the least favorable
submodel. The score functions for~$\beta$ and $q$ in the reparametrized
model are denoted by $\dot{\ell}_{1}(s,x;\beta,q)\,{=}\,(\partial/\partial
\beta) \log p^*_s(x;\beta,q) $ and $\dot{\ell}_{2}(s,x;\beta,q)=
(\partial/ \partial q) \log p^*_s(x;\beta,q)$, respectively.\vspace*{-2pt}

\begin{rem}
 In general, we may not have
the condition
\[
\int p^*_s(x;\beta,q)\,\mathrm{d}x =1\qquad
\mbox{for all } (\beta,q) \in \Theta_{\beta} \times D_q\ (s=1,\ldots,S).
\]
Therefore, there is no guarantee that each $p^*_s(x;\beta,q)$ is a
probability model. However,~(\ref{EqnDefreparameterized2}) ensures that
the linear combination $\sum_{s=1}^S w_s p^*_s(x;\beta,q)$ acts like a
probability model. This looks like a mixture model. The main
differences between the multisample model and the mixture model
are\vadjust{\goodbreak}
data and asymptotics. For example, the log-likelihood and the
information matrix in the mixture model are, respectively,
\[
\sum_{i=1}^n \log \Biggl\{ \sum_{s=1}^s w_s
p_s(x_i;\beta,q)\Biggr\}\vspace*{-2pt}
\]
and
\[
\int \biggl(\frac{(\partial/\partial (\beta,q))\sum_{s=1}^Sw_sp_s(x;\beta,q)}
{\sum_{s=1}^Sw_sp_s(x;\beta,q)}\biggr)^{\otimes
2}\sum_{s}w_sp_s(x;\beta,q)\,\mathrm{d}x,\vadjust{\goodbreak}
\]
while the log-likelihood and the information matrix in the multisample
model are given by, respectively, (\ref{log-likelihood_multi}) and
\[
\sum_{s=1}^Sw_s \int \biggl(\frac{(\partial/\partial
(\beta,q))p_s(x;\beta,q)}{p_s(x;\beta,q)}\biggr)^{\otimes 2}p_s(x;\beta,q)\,\mathrm{d}x.
\]
\end{rem}

\begin{rem}
Note that, since $q_{\beta_0}=\hat{\eta}_{\beta_0}=\eta_0$, we have
$p_s(x;\beta_0,\eta_0)= p^*_s(x;\beta_0,q_{\beta_0})$ ($s=1,\ldots,S$).
Therefore, for the reparametrized model, the notation $E_{s,0}$,
$s=1,\ldots,S$ is used for the expectations at the true value
$(\beta_0,q_{\beta_0})$.
\end{rem}

For a measurable function $f(s,x;\beta,q)$, define the
\textit{centering} of $f(s,x;\beta,q)$ by
\[
f^c(s,x;\beta,q)=f(s,x;\beta,q)-E_{s,0}\{f(s,x;\beta_0,q_{\beta_0})\}.
\]
The function $f^c(s,x;\beta,q)$ is called the \textit{centered}
$f(s,x;\beta,q)$.

\begin{theorem}[(Efficiency in a reparametrized model)]\label{thm31}
We assume that the least favorable submodel and the corresponding
reparametrized model are as in (\ref{ConditionR1}),
(\ref{ConditionR2}), (\ref{EqnDefreparameterized1}) and
(\ref{EqnDefreparameterized2}). Further, assume that
%
\begin{equation}\label{ConditionR3}
\frac{\partial}{\partial q}\bigg|_{q=q_{\beta}}\sum_{s=1}^S w_s
E_{s,0}\{\log p^*_s(x;\beta,q)\}=0 \qquad\mbox{for }  \beta \in
\Theta_{\beta}
\end{equation}
and
$\sum_{s=1}^S w_s E_{s,0}(\dot{\ell}_{2}^c \dot{\ell}_{2}^{cT})$
is non-singular. Then the efficient score function and the efficient
information matrix in the original multisample model $(\mathcal{P}_1,\ldots,\mathcal{P}_s)$
 are given by
\begin{equation}\label{reparametefficientscore}
  \dot{\ell}^{*}(s,x)
 = \dot{\ell}_{1}^c-
 \Biggl\{\sum_{s=1}^S w_s E_{s,0}(\dot{\ell}_{1}^c \dot{\ell}_{2}^{cT})\Biggr\}
 \Biggl\{\sum_{s=1}^S w_s E_{s,0}(\dot{\ell}_{2}^c \dot{\ell}_{2}^{cT})\Biggr\}^{-1}\dot{\ell}_{2}^c
\end{equation}
and
%
\begin{eqnarray}\label{reparametefficientinformation}
 I^{*} & = & \sum_{s=1}^S w_s E_{s,0}(\dot{\ell}_{1}^c \dot{\ell}_{1}^{cT}) \nonumber\\ [-8pt]\\ [-8pt]
& &{}  -  \Biggl\{\sum_{s=1}^S w_s E_{s,0}(\dot{\ell}_{1}^c
\dot{\ell}_{2}^{cT})\Biggr\} \Biggl\{\sum_{s=1}^S w_s
E_{s,0}(\dot{\ell}_{2}^c \dot{\ell}_{2}^{cT})\Biggr\}^{-1}
\Biggl\{\sum_{s=1}^S w_s E_{s,0}(\dot{\ell}_{2}^c
\dot{\ell}_{1}^{cT})\Biggr\},\nonumber\vadjust{\goodbreak}
\end{eqnarray}
where $\dot{\ell}_{1}^c(s,x;\beta,q)$ and
$\dot{\ell}_{2}^c(s,x;\beta,q)$ are the centered score functions for
$\beta$ and $q$ in the reparametrized model, respectively.
\end{theorem}

\begin{pf}
By (\ref{ConditionR2}) and (\ref{EqnDefreparameterized1}), the
efficient score function is given by
\begin{equation}\label{EqnTheoSpecialForm3}
\dot{\ell}^*(s,x)  =  \frac{\partial}{\partial
\beta}\bigg|_{\beta=\beta_0} \log p^*_s(x;\beta,q_{\beta})
 =  \dot{\ell}_{1}(s,x;\beta_0,q_{\beta_0})
+ \dot{q}_{\beta_0}^T\dot{\ell}_{2}(s,x;\beta_0,q_{\beta_0}).
\end{equation}

Since $E_{s,\beta_0\eta_0}\{\dot{\ell}^*(s,X)\}=0$ ($s=1,\ldots,S$), we
have
%
\begin{equation}\label{EqnTheoSpecialForm4}
E_{s,\beta_0\eta_0}\{\dot{\ell}_{1}(s,x;\beta_0,q_{\beta_0})\} +
\dot{q}_{\beta_0}^T
E_{s,\beta_0\eta_0}\{\dot{\ell}_{2}(s,x;\beta_0,q_{\beta_0})\}=0\qquad
(s=1,\ldots,S).
\end{equation}
Therefore, (\ref{EqnTheoSpecialForm3}) and (\ref{EqnTheoSpecialForm4})
imply
\begin{equation}\label{EqnTheoSpecialForm6}
\dot{\ell}^*(s,x) =  \dot{\ell}_{1}^c(s,x;\beta_0,q_{\beta_0}) +
\dot{q}_{\beta_0}^T\dot{\ell}_{2}^c(s,x;\beta_0,q_{\beta_0}).
\end{equation}

By differentiating (\ref{EqnDefreparameterized2}) with respect to $q$,
for all $(\beta,q) \in \Theta_{\beta} \times D_q$, we have
\[
\sum_{s=1}^S w_s \int \dot{\ell}_{2}(s,x;\beta,q) p^*_s(x;\beta,q)\,\mathrm{d}x=0.
\]
In particular, for all $\beta \in \Theta_{\beta}$,
\[
\sum_{s=1}^S w_s \int \dot{\ell}_{2}(s,x;\beta,q_{\beta}) p^*_s(x;\beta,q_{\beta})\,\mathrm{d}x=0.
\]
By differentiating with respect to $\beta$ at $\beta_0$,
\begin{eqnarray*}
& & \sum_{s=1}^S w_s \int
\biggl(\frac{\partial}{\partial \beta}\bigg|_{\beta=\beta_0}\dot{\ell}_{2}(s,x;\beta,q_{\beta})\biggr) p^*_s(x;\beta_0,q_{\beta_0})\,\mathrm{d}x\\
&&\quad =  - \sum_{s=1}^S w_s \int \dot{\ell}_{2}(s,x;\beta_0,q_{\beta_0})
\biggl(\frac{\partial}{\partial
\beta}\bigg|_{\beta=\beta_0}p^*_s(x;\beta,q_{\beta})\biggr)\,\mathrm{d}x.
\end{eqnarray*}
%
By the first equality in (\ref{EqnTheoSpecialForm3}), this equation is
equivalent to
%
\begin{equation}\label{EqnTheoSpecialForm5}
\sum_{s=1}^S w_s E_{s,0}\biggl\{\frac{\partial}{\partial
\beta}\bigg|_{\beta=\beta_0}\dot{\ell}_{2}(s,x;\beta,q_{\beta})\biggr\}
=- \sum_{s=1}^S w_s E_{s,0}(\dot{\ell}_{2}\dot{\ell}^{*T}).
\end{equation}

By differentiating (\ref{ConditionR3}) with respect to $\beta$ at
$\beta_0$, we get
\begin{eqnarray*}
0 & = & \frac{\partial}{\partial \beta}\bigg|_{\beta=\beta_0}
\frac{\partial}{\partial q}\bigg|_{q=q_{\beta}}\sum_{s=1}^S w_s
E_{s,0}\{\log p^*_s(x;\beta,q)\}
 =   \sum_{s=1}^S w_s E_{s,0}\biggl\{\frac{\partial}{\partial \beta}\bigg|_{\beta=\beta_0} \dot{\ell}_2(s,x,\beta,q_{\beta})\biggr\}\\
& = &  - \sum_{s=1}^S w_s E_{s,0}(\dot{\ell}_2\dot{\ell}^{*T})
 =   - \sum_{s=1}^S w_s E_{s,0}( \dot{\ell}_2^c\dot{\ell}^{*T}),
\end{eqnarray*}
where we used (\ref{EqnTheoSpecialForm5}) and
$E_{s,0}\{\dot{\ell}^{*}(s,X)\}=0$ ($s=1,\ldots,S$).

Therefore, the centered score function
$\dot{\ell}_2^c(s,x;\beta_0,q_{\beta_0})$ and the efficient score
function $\dot{\ell}^{*}(s,x)$ are uncorrelated.
Since $\dot{\ell}^{*} = \dot{\ell}_{1}^c +
\dot{q}_{\beta_0}^{T}\dot{\ell}_{2}^c$ (cf.
(\ref{EqnTheoSpecialForm6})), by the projection theorem (Theorem \ref{unprojectionmulti} in
Appendix \ref{appendixa}), we have
%
\begin{eqnarray*}
\dot{q}_{\beta_0}^T\dot{\ell}_{2}^c
& = & - \Biggl\{\sum_{s=1}^S w_s E_{s,0}(\dot{\ell}_{1}^c
\dot{\ell}_{2}^{cT})\Biggr\} \Biggl\{\sum_{s=1}^S w_s
E_{s,0}(\dot{\ell}_{2}^c\dot{\ell}_{2}^{cT})\Biggr\}^{-1}\dot{\ell}_{2}^c.
\end{eqnarray*}
The rest of the claims  follow by substituting this expression into
(\ref{EqnTheoSpecialForm6}).
\end{pf}

\begin{rem}
Under the usual regularity conditions, the solution
$(\hat{\beta}_n,\hat{q}_n)$ to the system of the score equations,
\[
\cases{\displaystyle
\sum_{s=1}^S\sum_{i=1}^{n_i}\dot{\ell}_{1}(s,X_{si};\hat{\beta}_n,\hat{q}_n)  =
0,\vspace*{2pt}\cr
\displaystyle\sum_{s=1}^S\sum_{i=1}^{n_i}\dot{\ell}_{2}(s,X_{si};\hat{\beta}_n,\hat{q}_n)
 =  0,
}
\]
is asymptotically distributed as
\[
\left\{\matrix{n^{1/2}(\hat{\beta}_n -\beta_0) \vspace*{2pt}\cr
n^{1/2}(\hat{q}_n-q_0)}\right\} \sim N\biggl\{\left(\matrix{0 \vspace*{2pt}\cr 0}\right), \Sigma^{-1}\biggr\},
\]
where
\[
\left.\Sigma=\cases{\displaystyle
\sum_{s=1}^Sw_s E_{s,0}(\dot{\ell}_{1}^c\dot{\ell}_{1}^{cT}), \sum_{s=1}^Sw_s E_{s,0}(\dot{\ell}_{1}^c\dot{\ell}_{2}^{cT})
\vspace*{2pt}\cr
\displaystyle\sum_{s=1}^Sw_s E_{s,0}(\dot{\ell}_{2}^c\dot{\ell}_{1}^{cT}),
\sum_{s=1}^Sw_s E_{s,0}(\dot{\ell}_{2}^c\dot{\ell}_{2}^{cT})
}\right\}.
\]
Then the asymptotic variance of
$n^{1/2}(\hat{\beta}_n -\beta_0)$ is given by $(I^{*})^{-1}$, where
$I^{*}$ is the efficient information for $\beta$ given by
(\ref{reparametefficientinformation}) (cf. Bickel \textit{et al.} \citet{BKRW1993}, page 28).
In this case, the estimator $\hat{\beta}_n$ is efficient. This
efficiency of the estimator based on the reparametrization is
demonstrated in a numerical example given in Section \ref{sec4}.
\end{rem}

\subsection{Example: Stratified sampling (continued)}\label{sec31}

In this section, we illustrate the use of Theorem \ref{thm31} to derive the
expressions of the efficient score function and the efficient
information bound in terms of the parameters in a reparametrized form
of the least favorable submodel in the stratified sampling
example.\looseness=1

Lemma \ref{lem21} gives the least favorable submodel with densities
\[
p_s(y,x;\theta,\hat{g}_{\theta})
 =   \frac{f(y|x;\theta)1_{(y,s) \in \mathcal{S}_s} \hat{g}_{\theta}(x)}{\hat{Q}_{s}(\theta)}\qquad  (s=1,\ldots,S),
\]
where $\hat{g}_{\theta}$ is given by (\ref{EqnhatgstratifiedSampling}).
By replacing
$\hat{Q}(\theta)=(\hat{Q}_1(\theta),\ldots,\hat{Q}_{S-1}(\theta),\hat{Q}_{S}(\theta))$
with $q=(q_1,\ldots,q_{S-1},1)$, we consider a reparametrized model of
the form
%
\begin{equation}\label{EqnReparameterizedModelStratified}
p^*_s(y,x;\theta,q)
 =   \frac{f(y|x;\theta)1_{(y,s) \in \mathcal{S}_s} \hat{g}_{\theta,q}(x)}{q_s}\qquad  (s=1,\ldots,S),
\end{equation}
where
%
\begin{equation}\label{EqnProofEffHatg}
\hat{g}_{\theta,q}(x)
  =  \frac{f^*_0(x)}
{ \sum_{s=1}^S  w_s Q_{s|X}(x;\theta)/q_s  }
\end{equation}
with $f^*_0(x)$ given by (\ref{Eqnf_star}).

The true value of $(\theta,q)$ is
\[
(\theta_0,q_0)=\biggl(\theta_0,\biggl(\frac{Q_1(\theta_0,g_0)}{Q_S(\theta_0,g_0)},
\ldots,\frac{Q_{S-1}(\theta_0,g_0)}{Q_S(\theta_0,g_0)},1\biggr)\biggr).
\]
Let $D_q$ be some neighborhood of $q_0$.

We will demonstrate that the conditions in Theorem \ref{thm31} are satisfied,
so that we can apply the theorem to identify the efficient score
function and the efficient information matrix in the example.

First,  we will show that
\[
\sum_{s=1}^S w_s \int p^*_s(y,x;\theta,q)\,\mathrm{d}y\,\mathrm{d}x=1\qquad
\mbox{for all } (\theta,q) \in \Theta_0 \times D_q.
\]

For any $(\theta,q)$, since  $Q_{s|X}(x;\theta)=\int
f(y|x;\theta)1_{(y,s) \in \mathcal{S}_s}\,\mathrm{d}y$,
\begin{eqnarray*}
\sum_{s=1}^S w_s \int p^*_s(y,x;\theta,q)\,\mathrm{d}y\,\mathrm{d}x
& = &\sum_{s=1}^S w_s \int \frac{f(y|x;\theta)1_{(y,s) \in \mathcal{S}_s} \hat{g}_{\theta,q}(x)}{q_s}\,\mathrm{d}y\,\mathrm{d}x\\
& = &\sum_{s=1}^S w_s \int \frac{Q_{s|X}(x;\theta)\hat{g}_{\theta,q}(x)}{q_s}\,\mathrm{d}x\\
& = & \int \sum_{s=1}^S w_s \frac{Q_{s|X}(x;\theta)}{q_s} \hat{g}_{\theta,q}(x)\,\mathrm{d}x\\
& = & \int f^*_0(x)\,\mathrm{d}x\qquad (\mbox{by (\ref{EqnProofEffHatg})})\\
& = & 1.
\end{eqnarray*}

Second,  we will show that for all $\theta \in \Theta_0$,
%
\begin{equation}\label{EqnCondEff}
\frac{\partial}{\partial q}\bigg|_{q=\hat{Q}(\theta)} \sum_{s=1}^Sw_s
E_{s,0}\{\log p_s(y,x;\theta,q)\}=0.
\end{equation}
For $j=1,\ldots,S-1$, the derivative is
\begin{eqnarray*}
& & \frac{\partial}{\partial q_j} \sum_{s=1}^Sw_s E_{s,0}\{\log p_s(y,x;\theta,q)\}\\
&&\quad =  - \frac{\partial}{\partial q_j} \sum_{s=1}^Sw_s
E_{s,0}\Biggl\{\log \sum_{s'=1}^S  w_{s'}
\frac{Q_{s'|X}(x;\theta)}{q_{s'}}
 +  \log  q_s \Biggr\}\\
&&\quad =  \sum_{s=1}^Sw_s
E_{s,0}\biggl\{\frac{w_jQ_{j|X}(x;\theta)/q_j^2}{ \sum_{s'=1}^S
w_{s'}  Q_{s'|X}(x;\theta)/q_{s'}} \biggr\} -  \frac{w_j}{q_j} \\
&&\quad =   \sum_{s=1}^S w_s \int \frac{w_jQ_{j|X}(x;\theta)/q_j^2}{
\sum_{s'=1}^S  w_{s'} Q_{s'|X}(x;\theta)/q_{s'}}
\frac{Q_{s|X}(x;\theta_0)g_0(x)}{Q_s(\theta_0,g_0)}\,\mathrm{d}x- \frac{w_j}{q_j} \\
&&\quad =   \int \frac{w_jQ_{j|X}(x;\theta)/q_j^2f^*_0(x)}{
\sum_{s'=1}^S  w_{s'} Q_{s'|X}(x;\theta)/q_{s'}}\,\mathrm{d}x
  -  \frac{w_j}{q_j}\qquad  (\mbox{by (\ref{Eqnf_star})})\\
&&\quad =  \frac{w_j}{q_j^2} \biggl(\int
Q_{j|X}(x;\theta)\hat{g}_{\theta,q}(x)\,\mathrm{d}x- q_j \biggr).
\end{eqnarray*}
Therefore, at
$q=(q_1,\ldots,q_{S-1},1)=(\frac{\hat{Q}_1(\theta)}{\hat{Q}_S(\theta)},\ldots,\frac{\hat{Q}_{S-1}(\theta)}{\hat{Q}_S(\theta)},1)$,
we have (\ref{EqnCondEff}).\vspace*{2pt}


By Theorem \ref{thm31},
 the efficient score function and the efficient information matrix in the example
 are calculated by (\ref{reparametefficientscore}) and (\ref{reparametefficientinformation}), respectively, where
the score functions are given by
\[
\dot{\ell}_{1}(s,y,x;\theta,q)= \frac{ (\partial/\partial
\theta)f(y|x;\theta)}{f(y|x;\theta)}
 - \frac{\sum_{s'=1}^S  w_{s'}
(\partial/\partial \theta)Q_{s'|X}(x;\theta)/q_{s'}}{ \sum_{s'=1}^S  w_{s'} Q_{s'|X}(x;\theta)/q_{s'}} \\
\]
and $\dot{\ell}_{2}(s,y,x;\theta,q)=
\{\dot{\ell}_{21}(s,y,x;\theta,q),\ldots,
\dot{\ell}_{2(S-1)}(s,y,x;\theta,q)\},$ where
\[
\dot{\ell}_{2j}(s,y,x;\theta,q)  =
  \frac{w_j}{q_{j}^2} \biggl\{\frac{Q_{j|X}(x;\theta)}{ \sum_{s'=1}^S  w_{s'}
  Q_{s'|X}(x;\theta)/q_{s'}} -q_{j}\biggr\}\qquad   (j=1,\ldots,S-1).
\]

Here verification of the non-singularity of
$\sum_{s=1}^Sw_sE_{s,0}(\dot{\ell}_{2}^c\dot{\ell}_{2}^{cT})$ is
omitted.

\section{Numerical example: Stratified sampling with logistic~regression}\label{sec4}


Here we compare the maximum likelihood estimator (MLE) and estimators
based on reparametri\-zations of the least favorable submodel, and
demonstrate that the estimators based on reparametri\-zations are
statistically as efficient as the MLE and computationally more
efficient.

The data in the Table \ref{tab1} were taken from Scott and Wild \citet{ScottWild1997,ScottWild2001} and
were the case-control sampling part of the study of people under 35 in
Northern Malawi.
Cases are those with new cases of leprosy and controls are those
without leprosy. The variable ``Scar'' indicates the presence or
absence of a BCG vaccination scar ($1 ={}$present, $0 ={}$absent).\

\begin{table}
\tablewidth=260pt
\caption{Leprosy data}\label{tab1}
\begin{tabular*}{260pt}{@{\extracolsep{\fill}}lllllll@{}}
\hline
&  \multicolumn{2}{l}{Scar${}=0$} & \multicolumn{2}{l}{Scar${}=1$} &
\multicolumn{2}{l@{}}{Total}\\ [-7pt]
&  \multicolumn{2}{l}{\hrulefill} & \multicolumn{2}{l}{\hrulefill} &
\multicolumn{2}{l@{}}{\hrulefill}\\
Age & Case & Control & Case & Control & Case & Control \\
\hline
\phantom{0}2.5  & \phantom{0}1  & 24 & \phantom{0}1  & 31 & \phantom{0}2 & 55\\
\phantom{0}7.5  & 11 & 22 & 14 & 39 & 25 & 61\\
12.5 & 28 & 23 & 22 & 27 & 50 & 50\\
17.5 & 16 & \phantom{0}5  & 28 & 22 & 44 & 27\\
22.5 & 20 & \phantom{0}9  & 19 & 12 & 39 & 21\\
27.5 & 36 & 17 & 11 & \phantom{0}5  & 47 & 22\\
32.5 & 47 & 21 & \phantom{0}6  & \phantom{0}3  & 53 & 24\\ [5pt]
Total &   &    &    &    &  260 & 260\\
\hline
\end{tabular*}
\end{table}

Let $x=(x_1,x_2)$ with $x_1={\rm Scar}$ and $x_2=100({\rm
Age}+7.5)^{-2}$. We consider a stratified sampling (case-control
sampling) with the logistic regression model
%
\begin{equation}\label{Eqnlogisticmodel}
f(y|x;\alpha,\beta)=\frac{\exp\{y(\alpha+x^T
\beta)\}}{1+\exp(\alpha+x^T \beta)}\qquad   (y \in \{0,1\},x \in R^2)
\end{equation}
and the partition $\mathcal{Y}\times \mathcal{X} = (\{0\} \times \mathcal{X})
\cup (\{1\} \times \mathcal{X})$,
 where $\alpha \in R$ and $\beta \in R^2$.
In this case, with $s=0,1$,
\[
Q_s(\alpha, \beta,g)=\int f(y=s|x;\alpha, \beta)g(x)\,\mathrm{d}x
\]
and
\[
Q_{s|X}(x,\alpha,\beta)=f(y=s|x;\alpha,\beta).
\]

From (\ref{EqnReparameterizedModelStratified}) and
(\ref{EqnProofEffHatg}), a  reparametrized model for the multisample
model is
\begin{eqnarray*}
p^*_s(x;\alpha,\beta,\rho_1)
 & = &    \frac{(q_0/q_s)f(y=s|x;\theta)}{ \sum_{s'=0}^1  w_{s'}  (q_0/q_{s'})Q_{s'|X}(x;\alpha,\beta)  }f^*_0(x)\\
 & = &    \frac{\exp\{s(\alpha + \log\rho_1 +x^T \beta)\}}{ w_{0} +  w_{1} \exp\{(\alpha+\log\rho_1+x^T \beta)\}  }f^*_0(x),
\end{eqnarray*}
where $\rho_0=q_0/q_0=1$ and $\rho_1=q_0/q_1$. The parameters in the
model are not identifiable and the parameters $\alpha$ and $\rho_1$
cannot be estimated separately. By the proof in the stratified sampling
example in Section \ref{sec31}, the efficient information bound for
$(\alpha,\beta)$ is given by (\ref{reparametefficientinformation}) in
Theorem \ref{thm31} with
$\dot{\ell}_{1}(s,x;\alpha,\beta,\rho_1)=\{\partial/\partial
(\alpha,\beta)\}\log p^{*}_s(x;\alpha,\beta,\rho_1)$ and
$\dot{\ell}_{2}(s,x;\alpha,\beta,\rho_1)=\{\partial/\partial\rho_1\}\log
p^{*}_s(x;\alpha,\beta,\rho_1)$.\vadjust{\goodbreak} The estimator
$(\hat{\alpha},\hat{\beta},\hat{\rho}_1)$ based on this
non-identifiable reparametrization is the maximizer of the
log-likelihood
$\ell_n(\alpha,\beta,\rho_1)=\sum_{s=0}^1\sum_{i=1}^{n_s}\log
p^{*}_s(x_{si};\alpha,\beta,\rho_1)$.

To gain identifiability of the parameters, we let $\alpha^*=\alpha +
\log \rho_1$, and the model is further reparametrized as
\[
p^{*}_s(x;\alpha^*,\beta) =   \frac{\exp\{s(\alpha^* +x^T \beta)\}}{
w_{0} +  w_{1} \exp\{(\alpha^*+x^T \beta)\}  }f^*_0(x).
\]
If we treat the parameters $\alpha$ and $g$ in the original model as
nuisance parameters, Theorem \ref{thm31} gives the efficient information bound
for an estimator of the parameter $\beta$: it is
(\ref{reparametefficientinformation}) in Theorem \ref{thm31} with
$\dot{\ell}_{1}(s,x;\alpha^*,\beta)=(\partial/\partial \beta)\log
p^{*}_s(x;\alpha^*,\beta)$ and
$\dot{\ell}_{2}(s,x;\alpha^*,\beta)=(\partial/\partial \alpha^*)\log
p^{*}_s(x;\alpha^*,\beta)$. The proof\vspace*{-1pt} is similar to the one for the
stratified sampling example given above and, therefore, we omit it. The
estimator $(\hat{\alpha}^*,\hat{\beta})$ based on this identifiable
reparametrization is the maximizer of the log-likelihood for the data
$\ell_n(\alpha^*,\beta)=\sum_{s=0}^1\sum_{i=1}^{n_s}\log
p^{*}_s(x_{si};\alpha^*,\beta)$.

If $X$ takes values in $\{v_1,\ldots,v_K\}$,  let $g_k=g(v_k)$,
$k=1,\ldots,K$. Then the log-likelihood for a~single observation in the
model can be written as
\[
\log p_s(x;\alpha,\beta,g)  =  \log f(y=s|x;\alpha, \beta)+
\sum_{k=1}^K 1_{\{x=v_k\}}\log g_k - \log \sum_{k=1}^K
f(y=s|v_k;\alpha,\beta)g_k.
\]
The MLE $(\hat{\alpha},\hat{\beta},\hat{g}),$ where
$\hat{g}=(\hat{g}_1,\ldots,\hat{g}_K),$ is the maximizer of the
log-likelihood $\ell_n(\alpha,\beta,g)=\sum_{s=0}^1\sum_{i=1}^{n_s}\log
p_s(x_{si};\alpha,\beta,g)$.

\begin{table}[b]
\caption{Model fitting results for the leprosy data}\label{tab2}
\begin{tabular*}{\textwidth}{@{\extracolsep{\fill}}lllllll@{}}
\hline
 &  \multicolumn{2}{l}{\multirow{2}{40pt}[-4pt]{Maximum likelihood}} & \multicolumn{4}{l@{}}{Reparametrization} \\ [-7pt]
 &   && \multicolumn{4}{l@{}}{\hrulefill}\\
 &   && \multicolumn{2}{l}{Not identifiable} &
 \multicolumn{2}{l@{}}{Identifiable}\\ [-7pt]
&  \multicolumn{2}{l}{\hrulefill} & \multicolumn{2}{l}{\hrulefill} &
 \multicolumn{2}{l@{}}{\hrulefill}\\
 & Coef & SE & Coef & SE &  Coef & SE \\
\hline
Intercept   & \phantom{$-$}1.55720   & 94.52766  & \phantom{$-$}0.61334  & 8388784     & --         & --  \\
Age         & $-$0.30205   & \phantom{0}0.19737 & $-$0.30211  & \phantom{0}0.19737   & $-$0.30215  &  0.19736 \\
Scar        & $-$4.30992   & \phantom{0}0.57891 & $-$4.31017  & \phantom{0}0.57892   & $-$4.30988  &  0.57889
\\ [5pt]
Computation time & \hspace*{2pt}43.61 sec  &    &   \hspace*{7pt}2.80 sec  &    &  \hspace*{7.1pt}2.44 sec   \\
\hline
\end{tabular*}
\end{table}

For each case (non-identifiable reparametrization, identifiable
reparametrization and maximum likelihood), let $\theta_1$ be the
parameter of interest and $\theta_2$ be the nuisance parameter. Then an
estimated variance of the estimator (of the parameter of interest) is
given by the formula (\ref{reparametefficientinformation}) except that
each $\sum_s w_s E_{s,0} (\dot{\ell}^c_i \dot{\ell}_j^{cT})$
($i,j=1,2$) is replaced with the corresponding second-degree partial
derivative $-n^{-1}(\partial^2/\partial \theta_i\,\partial \theta_j^T)
\ell_n$.

Estimates of  regression coefficients and their standard error (SE) in these models are given in Table \ref{tab2}.
Note that in the maximum likelihood and non-identifiable
reparametrization, the intercept parameter is not identifiable.
Its estimates and the corresponding SE
are unreliable and unstable. Therefore, we do not look at estimates of the
intercept
parameter in these models. The estimated coefficients of ``Age'' and ``Scar''  and their SE are very similar
to each other among these models. This is consistent with the
prediction made by Theorem \ref{thm31} that reparametrization gives the
semi-parametric efficiency bound that is achieved by the MLE.

\begin{table}
\tablewidth=280pt
\caption{Relative efficiency with respect to the maximum
likelihood}\label{tab3}
\begin{tabular*}{280pt}{@{\extracolsep{\fill}}lll@{}}
\hline
 &   \multicolumn{2}{l@{}}{Reparametrization}\\ [-7pt]
&   \multicolumn{2}{l@{}}{\hrulefill}\\
 &  Not identifiable & Identifiable\\
 \hline
Age   & 0.99997 & 0.99992  \\
Scar  & 1.00005 & 0.99994   \\ [5pt]
Computation time  &   0.06421   &    0.05595  \\
\hline
\end{tabular*}
\end{table}

Table \ref{tab3} gives the relative efficiency of estimates in non-identifiable
reparametrization and identifiable reparametrization with respect to
the maximum likelihood, along with the relative efficiency in
computation times (which is defined as  the ratio of the corresponding
computation times). The table indicates that these reparametrizations
are statistically as efficient as, and computationally more efficient
than, the method of maximum likelihood.


\section{Discussion}

Theorem \ref{thm31} gives conditions under which the efficient score function
and the efficient information matrix can be expressed in terms of the
parameters in the reparametrized model, namely
(\ref{reparametefficientscore}) and
(\ref{reparametefficientinformation}), respectively. In Section \ref{sec4}, we
demonstrated that Theorem~\ref{thm31} can be used to show the efficiency of
estimators based on non-identifiable and identifiable
reparametrizations in the logistic regression model, and that these
estimators are computationally more efficient than the MLE.
The results of the paper can be used to find a reparametrization of the
least favorable submodel (or profile likelihood) that gives
statistically and computationally efficient estimators in multisample
models.

\begin{appendix}
\section{}\label{appendixa}

We define the Hilbert space, projection and the efficient score
function.

\subsection{Hilbert space and the projection}
Let $\mathcal{H}$ be the Hilbert space of $m$-dimensional measurable
functions with zero mean and finite variance:
\[
\mathcal{H}=\Biggl\{\psi(s,x)\dvt  E_{s,0}(\psi)=0\ (s=1,\ldots,S),
\sum_{s=1}^S w_s E_{s,0}(\psi^T\psi)<\infty \Biggr\}.\vadjust{\goodbreak}
\]
The covariance of $\psi, \phi \in \mathcal{H}$ is defined by $\textrm{cov}
(\psi, \phi) =\sum_{s=1}^S w_s E_{s,0}(\psi \phi^T)$. We say $\psi$
and~$\phi$ are uncorrelated if $\textrm{cov} (\psi,\phi)=0$. For a set of
functions $\mathcal{G}$ in $\mathcal{H}$, $\mathcal{G}^{\perp}$ is the set of
all functions $\psi \in \mathcal{H}$ with $\textrm{cov} (\psi,\phi)=0$ for
all $\phi \in \mathcal{G}$. The projection $\Pi(\psi|\mathcal{G})$ of $\psi
\in \mathcal{H}$ onto a~closed subspace $\mathcal{G}$ is characterized by
\[
\Pi(\psi|\mathcal{G}) \in \mathcal{G}\quad   \mbox{and}\quad
\psi - \Pi(\psi|\mathcal{G}) \in \mathcal{G}^{\perp}.
\]

For an arbitrary Banach space $\mathcal{B}$, let $\mathcal{B}^*$ be its dual.
Let $A\dvtx\mathcal{B}\rightarrow \mathcal{H}$ be a bounded linear operator and
$\psi \in \mathcal{H}$. The \textit{adjoint operator} $A^T\dvtx\mathcal{H}\rightarrow \mathcal{B}^*$ of $A\dvtx\mathcal{B}\rightarrow \mathcal{H}$ is
defined by the map
\[
(A^T\psi)(b)= \langle Ab, \psi \rangle=\sum_{s=1}^S w_s E_{s,0}\{(Ab)\psi^T\},\qquad b \in \mathcal{B}.
\]

Suppose that $(A^T A)^{-1}$ exists and let $\psi \in \mathcal{H}$. By the
projection theorem for an operator equation,
\[
\Pi(\psi|\overline{A(\mathcal{B})})=A(A^T A)^{-1} A^T \psi
\]
is
a projection of $\psi$ onto the closure $\overline{A(\mathcal{B})}$ of the
range of $A$.

\subsection{The projection theorem}

\begin{theorem}[(The projection theorem)]\label{unprojectionmulti}
Suppose $\phi(s,x)$ is an $l$-dimensional vector of measurable functions such that
\begin{longlist}
\item[(1)] for $s=1,\ldots,S$, $E_{s,0}(\phi)=0$;
\item[(2)] $\sum_{s=1}^S w_s E_{s,0}(\phi^T\phi)<\infty$;\vspace*{2pt}
\item[(3)] $\{\sum_{s=1}^S w_s E_{s,0}( \phi \phi^T )\}^{-1}$
exists.
\end{longlist}

Let $\mathcal{G}=\{A \phi\dvt A \in R^{m \times l}\}$ be the closed
subspace of $\mathcal{H}$ generated by $\phi$. Then, for each $\psi \in
\mathcal{H}$, the projection of $\psi$ onto the closed subspace $\mathcal{G}$
is given by
\[
\pi(\psi|\mathcal{G}) = \Biggl\{\sum_{s=1}^S w_s E_{s,0}( \psi \phi^T )\Biggr\} \Biggl\{\sum_{s=1}^S w_s E_{s,0}(
\phi \phi^T )\Biggr\}^{-1} \phi.
\]
\end{theorem}

\begin{pf}
The proof is similar to the one for the standard case.
\end{pf}

\subsection{The efficient score function}

Here, we give the definition of the efficient score function in a
multisample model.

We assume the log-likelihood function for a single observation
$\ell(s,x;\beta,\eta)$ (defined by (\ref{EqnLogLikelihoodOne})) is
continuously differentiable with respect to $\beta$ for all $\beta \in
\Theta_{\beta}$ and
 Hadamard differentiable with respect to $\eta$ for all $\eta \in \Theta_{\eta}$.
The \textit{score function} $\dot{\ell}(s,x;\beta,\eta)$ for $\beta$
and the \textit{score operator} $A(s,x;\beta,\eta)$ for $\eta$ in the
multisample model are the derivatives of the log-likelihood function
with respect to $\beta$ and $\eta$, respectively.

The tangent space for $\eta$ is the closure $\overline{A(\mathcal{B})}$ of
range of the score operator $A$ for $\eta$.

The uncorrelated complement of the score function $\dot{\ell}_{\beta}$
 with respect to the tangent space for $\eta$,
\[
\dot {\ell} ^* =\dot {\ell} - \Pi (\dot {\ell} |\overline{A(\mathcal{B})} ),
\]
is called the \textit{efficient score function} in the
 multisample model $(\mathcal{P}_1,\ldots,\mathcal{P}_S)$.

\subsection{Theorem to identify the efficient score function}

To verify that the function given by (\ref{ConditionR2}) is the
efficient score function, the following theorem may be useful.

\begin{theorem}\label{theorema2}
A path $t \rightarrow \eta_t$ is a continuously differentiable map in a
neighborhood of~$0$ such that $\eta_{t=0}=\eta_0$. Define
$\alpha_t=\eta_t-\eta_0$. If $\beta \rightarrow \hat{\eta}_{\beta}$ is
a differentiable function such that
%
\begin{equation}\label{eqnEffCond1}
\hat{\eta}_{\beta_0}=\eta_0
\end{equation}
and, for each $\beta \in \Theta_{\beta}$, and for each path $\eta_t$,
%
\begin{equation}\label{eqnEffCond2}
\frac{\partial}{\partial t}\bigg|_{t=0}\sum_{s=1}^S w_s
E_{s,0}\{\log p_s(x;\beta,\hat{\eta}_{\beta}+\alpha_t)\}=0,
\end{equation}
then the function
%
\begin{equation}\label{EffScore}
\dot{\ell}^*(s,x)=\frac{\partial}{\partial \beta}\bigg|_{\beta=\beta_0}
\log p_s(x;\beta,\hat{\eta}_{\beta})
\end{equation}
 is the efficient score function.
\end{theorem}

\begin{pf}
Condition (\ref{eqnEffCond2}) implies that
\begin{eqnarray}\label{eqnEffCond3}
0 & = & \frac{\partial}{\partial \beta}\bigg|_{\beta=\beta_0} \frac{\partial}{\partial t}\bigg|_{t=0}
\sum_{s=1}^S w_s E_{s,0}\{\log p_s(x;\beta,\hat{\eta}_{\beta}+\alpha_t)\}
\nonumber\\ [-8pt]\\ [-8pt]
& = &  \frac{\partial}{\partial t}\bigg|_{t=0} \sum_{s=1}^S w_s
E_{s,0}\biggl\{\frac{\partial}{\partial \beta}\bigg|_{\beta=\beta_0}\log
p_s(x;\beta,\hat{\eta}_{\beta}+\alpha_t)\biggr\}.\nonumber
\end{eqnarray}

By differentiating the identity
\[
\sum_{s=1}^S w_s\int \biggl\{\frac{\partial}{\partial \beta}
\log p_s(x;\beta,\hat{\eta}_{\beta}+\alpha_t)\biggr\}p_s(x;\beta,\hat{\eta}_{\beta}+\alpha_t)\,\mathrm{d}x=0
\]
with respect to $t$ at $t=0$ and $\beta=\beta_0$, we get
\begin{eqnarray}\label{eqnEffCond4}
\hspace*{-30pt}0& = & \frac{\partial}{\partial t}\bigg|_{t=0,\beta=\beta_0}\sum_{s=1}^S w_s
\int \biggl(\frac{\partial}{\partial \beta}\log p_s(x;\beta,\hat{\eta}_{\beta}+\alpha_t)\biggr)p(x;\beta,\hat{\eta}_{\beta}+\alpha_t)\,\mathrm{d}x \nonumber\\
\hspace*{-30pt}& = & \sum_{s=1}^S w_s E_{s,0} \biggl[\dot{\ell}^*(s,x) \biggl\{\frac{\partial}{\partial t}\bigg|_{t=0}\log p_s(x;\beta_0,\eta_t)\biggr\}\biggr]
\qquad(\mbox{we used (\ref{EffScore}) and}\nonumber\\
\hspace*{-30pt}&&\hphantom{\sum_{s=1}^S w_s E_{s,0} \biggl[\dot{\ell}^*(s,x)
\biggl\{\frac{\partial}{\partial t}\bigg|_{t=0}\log p_s(x;\beta_0,\eta_t)\biggr\}\biggr]
\qquad(} \mbox{$\hat{\eta}_{\beta_0}+\alpha_t=\eta_t$ by (\ref{eqnEffCond1})})\\
\hspace*{-30pt}&&{}    +\frac{\partial}{\partial t}\bigg|_{t=0} \sum_{s=1}^S w_s E_{s,0} \biggl\{\frac{\partial}{\partial \beta}
\bigg|_{\beta=\beta_0}\log p_s(x;\beta,\hat{\eta}_{\beta}+\alpha_t)\biggr\}  \nonumber\\
\hspace*{-30pt}& = & \sum_{s=1}^S w_s E_{s,0} \biggl[\dot{\ell}^*(s,x)
\biggl\{\frac{\partial}{\partial t}\bigg|_{t=0}\log
p_s(x;\beta_0,\eta_t)\biggr\}\biggr]\qquad   (\mbox{by
(\ref{eqnEffCond3})}).\nonumber
\end{eqnarray}
Let $c \in R^m$ be arbitrary. Then, it follows from (\ref{eqnEffCond4})
that the product $c'\dot{\ell}^*(s,x)$ is orthogonal to the nuisance
tangent space $\dot\mathcal{P}_{\eta}$, which is the closed linear span of
score functions of the form $\phi(s,x)=\frac{\partial}{\partial
t}\big|_{t=0}\log p_s(x;\beta_0,\eta_t)$. By (\ref{EffScore}) with
(\ref{eqnEffCond1}) , we have
\begin{eqnarray*}
\dot{\ell}^*(s,x)
& = & \frac{\partial}{\partial \beta}\bigg|_{\beta=\beta_0} \log
p_s(x;\beta,\eta_0)+
\frac{\partial}{\partial \beta}\bigg|_{\beta=\beta_0} \log p_s(x;\beta_0,\hat{\eta}_{\beta})\\
& = & \dot{\ell}_{\beta}(s,x)- \psi(s,x),
\end{eqnarray*}
where $\dot{\ell}_{\beta}(s,x)=\frac{\partial}{\partial
\beta}\big|_{\beta=\beta_0} \log p_s(x;\beta,\eta_0)$ and
$\psi(s,x)=-\frac{\partial}{\partial \beta}\big|_{\beta=\beta_0} \log
p_s(x;\beta_0,\hat{\eta}_{\beta})$. Finally,
$c'\dot{\ell}^*(s,x)=c'\dot{\ell}_{\beta}(s,x)-c'\psi(s,x)$ is
orthogonal to the nuisance tangent space $\dot\mathcal{P}_{\eta}$ and
$c'\psi(s,x) \in \dot\mathcal{P}_{\eta}$ implies that $c'\psi(s,x)$ is the
orthogonal projection of $c'\dot{\ell}_{\beta}(s,x)$ onto the nuisance
tangent space $\dot\mathcal{P}_{\eta}$. Since $c \in R^m$ is arbitrary,
the function $\dot{\ell}^*(s,x)$ given by (\ref{EffScore}) is the
efficient score function.
\end{pf}

\section{}\label{appendixb}

\subsection{Proof of Lemma 2.1}

\begin{pf}
We show that  $\sum_{s=1}^{S}w_s \int \log
p_s(y,x;\theta,\hat{g}_{\theta})\,\mathrm{d}F_{s0}$ satisfies conditions
(\ref{eqnEffCond1}) and~(\ref{eqnEffCond2}) in Theorem \ref{theorema2} in Appendix
\ref{appendixa} so that the claim follows from this theorem.

Condition (\ref{eqnEffCond1}) is verified in Remark \ref{rem22}.
Now we verify (\ref{eqnEffCond2}). Let $g_t(x)$ be a path in the space
of density functions with $g_{t=0}(x)=g_0(x)$. Define
$\alpha_t(x)=g_t(x)-g_0(x)$ and write\ $\alpha'_0(x)=(\mathrm{d}/\mathrm{d}t)|_{t=0}
\alpha_t(x)$. Then
\begin{eqnarray*}
&&  \frac{\partial}{\partial t}\bigg|_{t=0}\sum_{s=1}^S w_s
\int \log p_s(y,x;\theta,\hat{g}_{\theta}+\alpha_t)\,\mathrm{d}F_{s0}\\
&&\quad =  \frac{\partial}{\partial t}\bigg|_{t=0} \sum_{s=1}^S w_s
\biggl[\int\log \{\hat{g}_{\theta}(x)+\alpha_t(x)\}\,\mathrm{d}F_{s,0}  - \log Q_s(\theta,\hat{g}_{\theta}+\alpha_t)\biggr]\\
&&\quad =  \frac{\partial}{\partial t}\bigg|_{t=0} \Biggl[\int   \log
\{\hat{g}_{\theta}(x)+\alpha_t(x)\} f^*_0(x)\,\mathrm{d}x
  -  \sum_{s=1}^S w_s \log Q_s(\theta,\hat{g}_{\theta}+\alpha_t)\Biggr]\\
&&\quad =   \int   \frac{\alpha'_0(x)}{\hat{g}_{\theta}(x)} f^*_0(x)\,\mathrm{d}x  -
\sum_{s=1}^S  w_s \frac{\int Q_{s|X}(x;\theta)\alpha'_0(x)\,\mathrm{d}x}{
\hat{Q}_s(\theta)}=0
\end{eqnarray*}
by (\ref{EqnhatgstratifiedSampling}) and (\ref{Eqnf_star}).
\end{pf}
\end{appendix}


\printhistory

\end{document}